\documentclass[11pt]{article}

\usepackage{amssymb,amsmath}
\usepackage{latexsym}

{\catcode `\@=11 \global\let\AddToReset=\@addtoreset} \AddToReset{equation}
{section}

\newtheorem{lemma}{\bf Lemma}[section]
\newtheorem{corollary}{\bf Corollary}[section]
\newtheorem{proposition}{\bf Proposition}[section]

\newtheorem{@definition}{\sc Definition}[section]
\newtheorem{@remark}{\sc Remark}[section]
\newenvironment{remark}{\begin{@remark}\rm}{\end{@remark}}
\newtheorem{@example}{\sc Example}[section]
\newenvironment{example}{\begin{@example}\rm}{\end{@example}}
\def\mathsf{\bf}

\def\supp{\rm supp}

\def\E{\mathrm E}

\def\d{\mathrm d}
\def\e{\mathrm e}
\def\i{\mathrm i}

\def\text{\mbox}

\newcommand{\Cov}{{\rm Cov}}

\newcommand{\beq}{\begin{equation}}
\newcommand{\eeq}{\end{equation}}
\newcommand\beqn{\begin{displaymath}}  % no number
\newcommand\eeqn{\end{displaymath}}
\newcommand\refeq[1]{{\rm (\ref{e:#1})}}
\reversemarginpar

\begin{document}

\title{Time series aggregation, disaggregation and long memory
\footnote{The research is supported by joint Lithuania and France scientific program Gilibert, PAI EGIDE 09393 ZF.}}

\author{Dmitrij Celov${}^1$, Remigijus Leipus${}^{1,2}$
and Anne Philippe${}^3$  \\
{\small ${}^1$Vilnius University}\\
{\small ${}^2$Vilnius Institute of Mathematics and Informatics}\\
{\small ${}^3$Laboratoire de Math\'ematiques Jean Leray, Universit\'e de
  Nantes}
} \maketitle

\vskip2cm

\small{
\begin{quote}
{\bf  Abstract} The paper studies the aggregation/disaggregation problem of random parameter
 AR(1) processes and its relation to the long memory phenomenon.
 We give a characterization of a subclass of aggregated processes which can be obtained
 from simpler, "elementary", cases. In particular cases of the mixture densities, the structure
 (moving average representation) of the aggregated process is investigated.
\end{quote}
}

\vskip.2cm

\vskip1cm

\vskip.5cm


\section{Introduction}

In this paper we consider the aggregation scheme introduced in the paper of Granger (1980), where it was shown
 that aggregation of random parameter AR(1) processes can produce long memory. Since this work, a large number
 of papers were devoted to the question how do micro level processes imply the long memory at macro level and
 applications (mainly in economics) (see Haubrich and Lo (2001), Oppenheim and Viano (2004), Zaffaroni (2004), Chong (2006), etc.).

Recently, Dacunha-Castelle and Oppenheim (2001) (see also
Dacunha-Castelle and Fermin (2006))
 stated the following problem: which class of long memory processes can be obtained by the aggregation
 of short memory models with random coefficients?

Let us give a precise formulation of the problem when the
underlying short memory models are described by AR(1) dynamics.
Let
 $\varepsilon=\{\varepsilon_t,t\in {\bf Z}\}$ be a
 sequence of independent identically distributed (i.i.d.) random variables (r.v.) with $\E \varepsilon_t=0$,
 $\E \varepsilon^2_t=\sigma^2_\varepsilon$, and let $a$ be a random variable supported by $(-1,1)$ and satisfying
\beq\label{e:as}
 \E \Big[\frac{1}{1-a^2}\Big]<\infty.
\eeq
 Consider a sequence of i.i.d.\ processes $Y^{(j)}=\{Y^{(j)}_t,t\in {\bf Z}\}$, $j\ge 1$
 defined by the random AR(1) dynamics
\beq\label{e:AR1}
 Y^{(j)}_t= a^{(j)} Y^{(j)}_{t-1} +\varepsilon_t^{(j)},
\eeq
 where $\varepsilon^{(j)}=\{\varepsilon^{(j)}_t,t\in {\bf Z}\}$, $j\ge 1$, are independent
 copies of $\varepsilon=\{\varepsilon_t,t\in {\bf Z}\}$ and $a^{(j)}$'s are independent copies of $a$.
 Here, the sequences $a, a^{(j)}$, $j\ge 1$ and $\varepsilon$, $\varepsilon^{(j)}$, $j\ge 1$, are independent.
 Under these conditions, \refeq{AR1} admits a covariance-stationary
solution $Y^{(j)}_t$
 and the finite dimensional distributions of the process $X^{(N)}_t=N^{-1/2}\sum_{j=1}^N Y^{(N)}_t$, $t\in {\bf Z}$,
 weakly converge, as $N\to\infty$, to those of a zero mean Gaussian process $X_t$, called the
 {\em aggregated process} (see Oppenheim and Viano (2004)).

Assume that the distribution of r.v.\ $a$ admits a {\em mixture
density} $\varphi$, which by \refeq{as}
 satisfies
\beq\label{e:disp}
 \int_{-1}^1 \frac{\varphi(x)}{1-x^2}\; \d x<\infty.
\eeq
 The covariance function and spectral density of aggregated process $X_t$ are given, respectively, by
\beq \label{e:cov}
 \gamma(h):=\Cov(X_h,X_0)=\sigma^2_\varepsilon
 \int^1_{-1} \frac{x^{|h|}}{1-x^2}\;\varphi(x) \d x
\eeq
 and
\beq\label{e:sp_density}
 f(\lambda)=\frac{\sigma^2_\varepsilon}{2\pi}\int^1_{-1}\frac{\varphi(x)}{|1-x\e^{\i\lambda}|^2}\; \d x.
\eeq
 Note that an aggregated process $X_t$ possess the long memory property
 (i.e.\ $\sum_{h=-\infty}^\infty |\gamma(h)|=\infty$) if and only if
\beq  \label{e:lm}
 \int_{-1}^1 \frac{\varphi(x)}{(1-x^2)^2} \;\d x =\infty.
\eeq

If the mixture density $\varphi$ is {\em a priori given} and our
aim is to characterize the properties of the induced aggregated
process (moving average representation, spectral density,
covariance function, etc.), we
 call this problem an {\em aggregation problem}. And vice versa, if we observe the aggregated process
 $X_t$ with spectral density $f$ and we need to find the individual processes (if they exist) of form \refeq{AR1} with some
 mixture density $\varphi$, which produce the aggregated process, then we call this problem a
 {\em disaggregation problem}. The second problem, which is much
  harder than the first one, is equivalent to the finding of $\varphi$ such that \refeq{sp_density} (or \refeq{cov}) and \refeq{disp} hold.
  In the latter case we say that {\em the mixture density $\varphi$
 is associated with the spectral density $f$}.

Sections \ref{s:product} and \ref{s:slm} are devoted to the disaggregation problem. Equality \refeq{cov} shows that the covariance function
$\gamma(h)$ can be interpreted as a $h$-moment of the density function $\varphi(x)(1-x^2)^{-1}$ supported by $(-1,1)$, and thus finding of the
mixture density is related to the moments' problem (see Feller (1976)). In Section~\ref{s:product}, we prove that, under some mild conditions,
the mixture density associated with the product spectral density can be obtained from the "elementary" mixture densities associated with the
multipliers. In Section~\ref{s:slm}, we apply the obtained result to the spectral density which is a product of to FARIMA-type spectral
densities. Using the form of mixture density for the classical FARIMA model we show that one can solve the disaggregation problem for more
complex long memory stationary sequences.

For the aggregation problem, it is important to characterize, for
a given class of mixture densities, the behavior of the
coefficients in its linear representation and the behavior of the
spectral density. We address this problem to Section~\ref{s:MA}.
In the Appendix, we provide a proof of the form of the mixture
density in FARIMA case.

\section{Mixture density for the product of aggregated spectral densities}
\label{s:product}

%{\bf [[Assume that the spectral densities in \refeq{ff12} are of
%the forms $f_1(\lambda)=(2\pi)^{-1} |F_1(\e^{i\lambda})|^{-2}$,
%$f_2(\lambda)=(2\pi)^{-1} |F_2(\e^{i\lambda})|^{-2}$ and the
%processes $X^{(1)}_t$ and $X^{(2)}_t$, defined as
%$F_j(B)X^{(j)}_t=\varepsilon^{(j)}_t$, $j=1,2$, are obtained by
%aggregating the AR(1) models with mixture densities $\varphi_1$
%and $\varphi_2$, respectively. Proposition~\ref{p:f12} implies
%that then the process $X_t$ given by
%$F_1(B)F_2(B)X_t=\varepsilon_t$, where $B$ is backshift operator,
%is aggregated process as well, and it's mixture density is given
%in \refeq{mix}.]]}

%In this section we show how the mixture density associated with
% the product of spectral densities can be obtained from (simpler) mixture
% densities associated with multipliers.

 Let $X_{1,t}$ and $X_{2,t}$
 are two aggregated processes obtained from the independent copies of
 AR(1) sequences $Y_{1,t}= a_1 Y_{1,t-1} +\varepsilon_{1,t}$ and
 $Y_{2,t}= a_2 Y_{2,t-1} +\varepsilon_{2,t}$, respectively, where $a_1$, $a_2$
 satisfy $\refeq{as}$, $\varepsilon_{i}=\{\varepsilon_{i,t},t\in {\bf Z}\}$
 and $a_i$ are independent, $i=1,2$. Denote $\sigma_{i,\varepsilon}^2:=\E\varepsilon^2_{i,t}$, $i=1,2$.

 Assume that $\varphi_1$ and $\varphi_2$ are the mixture densities associated
 with spectral densities $f_1$ and $f_2$, respectively, i.e.
\beq\label{e:aaa}
 f_i(\lambda)=\frac{\sigma^2_{i,\varepsilon}}{2\pi}\int^1_{-1}
 \frac{\varphi_i(x)}{|1-x\e^{\i\lambda}|^2}\; \d x, \ \ i=1,2.
\eeq
 The following proposition shows that, if $\varphi_1$ and $\varphi_2$ in \refeq{aaa}
 are supported by $[0,1]$ and $[-1,0]$ respectively, then
  the stationary Gaussian
 process with spectral density $f(\lambda)=f_1(\lambda)f_2(\lambda)$
 {\em can also be obtained by aggregation of the i.i.d.\ AR(1) processes}
 with some mixture density $\varphi$ and noise sequence
 $\varepsilon$.

\begin{proposition}\label{p:f12}
Let $\varphi_1$ and $\varphi_2$ be the mixture densities
associated
 with spectral densities $f_1$ and $f_2$, respectively. Assume that
$\supp(\varphi_1) \subset [0,1]$, $\supp(\varphi_2)
 \subset [-1,0]$ and
% Let $X_t, t\in {\bf Z}$ be a process having spectral density of the form
\beq\label{e:ff12}
 f(\lambda) = f_1(\lambda) f_2(\lambda).
\eeq
   Then the mixture density $\varphi(x)$, $x\in [-1,1]$ associated with $f$ is given by equality
\beq\label{e:mix}
 \varphi(x) =\frac{1}{ C_*}\bigg(\varphi_1(x)
 \int_{-1}^{0}\frac{\varphi_2(y)}{(1-xy)(1-y/x)} \;\d y +
 \varphi_2(x) \int_{0}^{1}\frac{\varphi_1(y)}{(1-xy) (1-y/x)}  \;\d
 y\bigg),
\eeq
 where $C_*:=\int_0^1\Big(\int_{-1}^0
 \varphi_1(x)\varphi_2(y)(1-xy)^{-1}\d y\Big)\d x$. The
 variance of the noise is
\beq\label{e:sigma2epsilon}
  \sigma^2_\varepsilon=\frac{\sigma^2_{1,\varepsilon} \sigma^2_{2,\varepsilon} C_*}{2\pi}.
\eeq
\end{proposition}

%\begin{eqnarray*}
%&&\int_{-1}^1\bigg( \varphi_1(x)
% \int_{-1}^{0}\frac{\varphi_2(y)}{(1-xy)(1-y/x)} \;\d y +
% \varphi_2(x) \int_{0}^{1}\frac{\varphi_1(y)}{(1-xy) (1-y/x)}  \;\d
% y\bigg)\d x\\
%&=&\int_0^1\bigg(\int_{-1}^{0}\frac{ \varphi_1(x)\varphi_2(y)}{(1-xy)(1-y/x)} \;\d y\bigg)\d x+
% \int_{-1}^0\bigg(\int_0^1\frac{\varphi_2(x)\varphi_1(y)}{(1-xy)(1-y/x)} \;\d y\bigg)\d x\\
%&=&\int_0^1\bigg(\int_{-1}^{0}\frac{ \varphi_1(x)\varphi_2(y)}{(1-xy)(1-y/x)} \;\d y\bigg)\d x+
% \int_0^1\bigg(\int_{-1}^{0}\frac{\varphi_2(y)\varphi_1(x)}{(1-xy)(1-x/y)} \;\d y\bigg)\d x\\
%&=&\int_0^1\bigg(\int_{-1}^{0}\frac{ \varphi_1(x)\varphi_2(y)}{1-xy} \bigg(\frac{1}{1-y/x}+\frac{1}{1-x/y}\bigg)\d y \bigg)\d x\\
%&=&\int_0^1\bigg(\int_{-1}^{0}\frac{\varphi_1(x)\varphi_2(y)}{1-xy}\; \d y \bigg)\d x=C_*.
%\end{eqnarray*}

%\newpage

\noindent {\sc Proof.} Obviously, the covariance function
$\gamma(h)=\Cov(X_h,X_0)$ has a form
%{\bf [[FORMALLY: we still
%don't know if this series converges]]}
\begin{equation}\label{eq:1}
    \gamma(h) = \frac{1}{2\pi}\sum_{j=-\infty}^\infty \gamma_1(j+|h|)
    \gamma_2(j),
\end{equation}
 where $\gamma_1$ and $\gamma_2$ are the covariance functions of
 the aggregated processes obtained from the mixture densities $\varphi_1$ and $\varphi_2$,
 respectively:
$$
 \gamma_1(j) =\sigma^2_{1,\varepsilon} \int_0^{1} \varphi_1(x)
  \frac{x^{|j|}}{1-x^2} \; \d x,\ \gamma_2(j) =\sigma^2_{2,\varepsilon} \int_{-1}^0 \varphi_2(x)
  \frac{x^{|j|}}{1-x^2} \; \d x.
$$
Clearly, $\gamma_1(j)>0$ and $\gamma_2(j)=(-1)^j |\gamma_2(j)|$.
Let $h\ge 0$. Then
\begin{eqnarray}
 2\pi\sigma^{-2}_{1,\varepsilon} \sigma^{-2}_{2,\varepsilon}\gamma(h) &=& \sigma^{-2}_{1,\varepsilon} \sigma^{-2}_{2,\varepsilon}
  \sum_{j=-\infty}^\infty \gamma_1(j+h) \gamma_2(j) \nonumber\\
  &=& \sigma^{-2}_{1,\varepsilon} \sigma^{-2}_{2,\varepsilon}\bigg(\sum_{j=0}^{\infty} \gamma_1(j+h) \gamma_2(j) + \sum_{j=0}^{\infty} \gamma_1(j)
  \gamma_2(j+h) + \sum_{j=1}^{h-1} \gamma_1(-j+h) \gamma_2(j)\bigg)\nonumber\\
  &=:& s_1+s_2+s_3.\label{e:s}
\end{eqnarray}
 We have $s_1=\lim_{N\to\infty} s_1^{(N)}$, where
\begin{align*}
  s_1^{(N)} &= \sum_{j=0}^{N} \int_{0}^{1} \varphi_1(x) \frac{x^{j+h}}{1-x^2}\;\d x
  \int_{-1}^{0} \varphi_2(y) \frac{y^{j}}{1-y^2}\; \d y \\
  &=  \int_{0}^{1}\int_{0}^{1}  \varphi_1(x)\varphi_2(-y)  \frac{1}{1-x^2}\frac{1}{1-y^2} \sum_{j=0}^{N} (-1)^j  x^{j+h} y^{j}\; \d x \;\d y\\
 &=  \int_{0}^{1}\int_{0}^{1}  \varphi_1(x)\varphi_2(-y)
 \frac{x^h}{1-x^2}\frac{1}{1-y^2}
\frac{1-(-xy)^{N+1}}{1+xy} \ \d x \;\d y.
\end{align*}
 Note that
 \beqn
  \left\vert \varphi_1(x)\varphi_2(-y) \frac{x^h}{1-x^2}\frac{1}{1-y^2}
    \frac{1-(-xy)^{N+1}}{1+xy} \right\vert \leq 2
  \varphi_1(x)\varphi_2(-y) \frac{1}{1-x^2}\frac{1}{1-y^2}
\eeqn
 and
\beqn
  \int_0^1\int_0^1 \varphi_1(x)\varphi_2(-y)
  \frac{1}{1-x^2}\frac{1}{1-y^2}\;\d x\;\d y = \gamma_1(0)\gamma_2(0) <\infty.
\eeqn
 Therefore, by the dominated convergence
theorem, we obtain
\begin{align}
  s_1&=
  \int_{0}^{1}\int_{0}^{1}   \varphi_1(x)\varphi_2(-y)
 \frac{x^h}{1-x^2}\frac{1}{1-y^2}
\lim_{N\to\infty}\frac{1-(-xy)^{N+1}}{1+xy}\; \d x \d y \nonumber \\
&= \int_{0}^{1}\int_{0}^{1}  \varphi_1(x)\varphi_2(-y)
 \frac{x^h}{1-x^2}\frac{1}{1-y^2}
\frac{1}{1+xy}\; \d x \d y \nonumber \\
&= \int_{0}^{1} \frac{x^h}{1-x^2}  \left\{  \varphi_1(x) \int_{-1}^{0}  \varphi_2(y) \frac{1}{1-y^2} \frac{1}{1-xy}\; \d y \right\} \d x.
\label{e:s1}
\end{align}

Analogously, we have
\beq
 s_2= \int_{-1}^{0}   \frac{y^h}{1-y^2} \left\{ \varphi_2(y)
 \int_{0}^{1}  \varphi_1(x) \frac{1}{1-x^2}
 \frac{1}{1-xy} \;\d x  \right\} \d y. \label{e:s2}
\eeq

For the last term in the decomposition \refeq{s} we have
\begin{eqnarray}
 s_3 &=&  \sum_{j=1}^{h-1}
\int_{0}^{1} \varphi_1(x) \frac{x^{-j+h}}{1-x^2}\; \d x \int_{-1}^{0} \varphi_2(x) \frac{x^{j}}{1-x^2}\; \d x\nonumber \\
&=&  \int_{0}^{1} \int_{0}^{1}  \varphi_1(x) \frac{x^{h}}{1-x^2}\;
\varphi_2(-y) \frac{1}{1-y^2} (-y/x) \frac{1 -(-y/x)^{h-1}}{1+y/x}\; \d x \;\d y\nonumber\\
 &=&  \int_{0}^{1}  \frac{x^{h}}{1-x^2} \left\{  \varphi_1(x)   \int_{-1}^{0}
\varphi_2(y) \frac{1}{1-y^2}  \frac{y/x}{1-y/x} \; \d y \right\} \d x \nonumber\\
&& -  \int_{-1}^{0} \frac{y^h}{1-y^2} \left\{ \varphi_2(y) \int_{0}^{1}  \varphi_1(x) \frac{1}{1-x^2}
  \frac{1}{1-y/x}\; \d x  \right\}\d y.\label{e:s3}
\end{eqnarray}

Equalities \refeq{s1}--\refeq{s3}, together with \refeq{s}, imply
\begin{eqnarray*}
2\pi\sigma^{-2}_{1,\varepsilon} \sigma^{-2}_{2,\varepsilon}
\gamma(h) &=& \int_{0}^{1} \frac{x^{h}}{1-x^2}
\;\varphi_1(x)\bigg\{
    \int_{-1}^{0}  \varphi_2(y) \frac{1}{1-y^2}  \frac{y/x}{1-y/x}\;  \d y \\
 && + \int_{-1}^{0}  \varphi_2(y)
\frac{1}{1-y^2}   \frac{1}{1-xy}\; \d y \bigg\}  \d x \\
&&+ \int_{-1}^{0} \frac{y^h}{1-y^2}\;\varphi_2(y) \bigg\{\int_{0}^{1}  \varphi_1(x) \frac{1}{1-x^2}
  \frac{1}{1-xy}\; \d x   \\
&& - \int_{0}^{1}  \varphi_1(x)
 \frac{1}{1-x^2}\frac{1}{1-y/x}\; \d x  \bigg\} \d y\\
&=& \int_{0}^{1} \frac{x^{h}}{1-x^2} \;\varphi_1(x)\bigg\{
    \int_{-1}^{0}  \varphi_2(y) \frac{1}{(1-y/x)(1-xy)}\;  \d y\bigg\}\d x \\
&&+ \int_{-1}^{0} \frac{y^h}{1-y^2}\;\varphi_2(y) \bigg\{\int_{0}^{1}  \varphi_1(x)
  \frac{1}{(1-x/y)(1-xy)}\; \d x  \bigg\}\d y.
\end{eqnarray*}
This and \refeq{cov} imply \refeq{mix}, taking into account that
$\int_{-1}^1\varphi(x)\d x=1$. \hfill $\Box$
\medskip

%{\bf [[Assume that the spectral densities in \refeq{ff12} are of
%the forms $f_1(\lambda)=(2\pi)^{-1} |F_1(\e^{i\lambda})|^{-2}$,
%$f_2(\lambda)=(2\pi)^{-1} |F_2(\e^{i\lambda})|^{-2}$ and the
%processes $X^{(1)}_t$ and $X^{(2)}_t$, defined as
%$F_j(B)X^{(j)}_t=\varepsilon^{(j)}_t$, $j=1,2$, are obtained by
%aggregating the AR(1) models with mixture densities $\varphi_1$
%and $\varphi_2$, respectively. Proposition~\ref{p:f12} implies
%that then the process $X_t$ given by
%$F_1(B)F_2(B)X_t=\varepsilon_t$, where $B$ is backshift operator,
%is aggregated process as well, and it's mixture density is given
%in \refeq{mix}.]]}

\section{Seasonal long memory case}\label{s:slm}

In this section, we apply the obtained result to the spectral
densities $f_1$ and $f_2$
%$$
% h(\lambda;d)=\frac{1}{2\pi}\ \Big(2\sin\frac{|\lambda|}{2}\Big)^{-2d}, \ \ 0<d<1/2.
%$$
 having the forms:
\begin{eqnarray}
 f_1(\lambda;d_1)&=& \frac{1}{2\pi}\ |1-\e^{\i \lambda}|^{-2d_1}% (2\sin\frac{|\lambda|}{2}\Big)^{-2d_1},
  \label{e:FARIMA}\\
 f_2(\lambda;d_2)&=& \frac{1}{2\pi}\ |1+\e^{\i \lambda}|^{-2d_2},
% \Big(2\cos\frac{\lambda}{2}\Big)^{-2d_2},
 \ \ 0<d_1,d_2<1/2.\label{e:SFARIMA}
\end{eqnarray}
 We call these spectral densities (and corresponding processes)
 {\em fractionally integrated}, FI($d_1$), and
 {\em seasonal fractionally integrated}, SFI($d_2$), spectral
 densities.
 % {\bf [[Maybe we shall call differently.]]}
 The
 mixture density associated with the FI($d_1$) spectral
 density \refeq{FARIMA} is given by the following
 expression (for the sketch of proof see Dacunha-Castelle and Oppenheim (2001)):
\beq\label{e:phi1}
 \varphi_1(x;d_1)=
 C(d_1) x^{d_1-1} (1-x)^{1-2d_1} (1+x) {\bf 1}_{[0,1]}(x)
\eeq
 with
\beq\label{e:constant:fd}
 C(d)=\frac{\Gamma(3-d)}{2\Gamma(d)\Gamma(2-2d)}
 = 2^{2d-2}\frac{\sin(\pi d)}{\sqrt \pi} \frac{\Gamma(3-d)}{\Gamma ((3/2)-d)}
\eeq
 and the variance of the noise
\beq\label{e:noise:fd}
  \sigma^2_{1,\varepsilon}= \frac{\sin(\pi d_1)} {C(d_1)\pi}.
\eeq
 For convenience, we provide the rigorous proof of this result in Proposition~\ref{p:DCO} of
 Appendix.

 Similarly, the mixture density associated with the
 spectral density \refeq{SFARIMA} is given by
\beq\label{e:phi2}
 \varphi_2(x;d_2)= \varphi_1(-x;d_2)=C(d_2) |x|^{d_2-1} (1+x)^{1-2d_2} (1-x) {\bf
 1}_{[-1,0]}(x)
\eeq
 and
\beq\label{e:snoise:fd}
  \sigma^2_{2,\varepsilon}= \frac{\sin(\pi d_2)} {C(d_2)\pi},
\eeq
 since
\beqn
 f_2(\lambda;d_2)=f_1(\pi-\lambda;d_2)= \frac{\sigma^2_{2,\varepsilon}} {2\pi} \int_0^1
 \frac{\varphi_1(x;d_2)}{|1+x\e^{\i\lambda}|^2}\;\d x=
 \frac{\sigma^2_{2,\varepsilon}} {2\pi} \int_{-1}^0
 \frac{\varphi_1(-x;d_2)}{|1-x\e^{\i\lambda}|^2}\;\d x.
\eeqn

Clearly, mixture densities $\varphi_1(x;d_1)$ and
$\varphi_2(x;d_2)$, given in
 \refeq{phi1}, \refeq{phi2}, satisfy \refeq{disp} and, hence, assumptions of
 Proposition~\ref{p:f12},
 whenever $0<d_1<1/2$ and $0<d_2<1/2$. Moreover, since $d_1>0$ and $d_2>0$, both $\varphi_1$ and
 $\varphi_2$ satisfy \refeq{lm}.
%and hence $f_1$ and $f_2$ are spectral densities of
%long memory processes.
The mixture density associated with the spectral density
\beq\label{e:ff}
 f(\lambda;d_1,d_2)=f(\lambda;d_1) f(\lambda;d_2)=\frac{1}{(2\pi)^2}\
 |1-\e^{\i \lambda}|^{-2d_1}  |1+\e^{\i \lambda}|^{-2d_2},\
% \frac{1}{(2\pi)^2}\
% \Big(2\sin\frac{|\lambda|}{2}\Big)^{-2d_1}
% \Big(2\cos\frac{\lambda}{2}\Big)^{-2d_2},
 0<d_1,d_2<1/2
\eeq
 can be derived from representation \refeq{mix}.

Denote $F(a,b,c;x)$ a hypergeometric function
$$
 F(a,b,c;x)=\frac{\Gamma(c)}{\Gamma(b)\Gamma(c-b)}\int_0^1
 t^{b-1}(1-t)^{c-b-1}(1-tx)^{-a} \d t,
$$
 where $c>b>0$ if $x<1$ and, in addition, $c-a-b>0$ if $x=1$.

\begin{proposition} %{\bf [[to check the constants in the proof]]}
 The mixture density associated with $f(\cdot;d_1,d_2)$ \refeq{ff} is given by equality
\begin{eqnarray}
 \varphi(x;d_1,d_2)&=&C(d_1,d_2) x^{d_1-1} (1-x)^{1-2d_1}
 G(-x;d_2){\bf 1}_{[0,1]}(x)\nonumber\\
 && + C(d_2,d_1) |x|^{d_2-1} (1+x)^{1-2d_2}
 G(x;d_1){\bf 1}_{[-1,0]}(x),\label{e:phif}
\end{eqnarray}
 where
$$
  G(x;d):= F\Big(1,d,2-d;\frac{1}{x}\Big)-x F(1,d,2-d;x)
$$
and
\begin{eqnarray*}
 C(d_1,d_2) &=& (C^*)^{-1} \frac{\Gamma(d_2)\Gamma(2-2d_2)}{\Gamma(2-d_2)},\\
 %= (C^*)^{-1}
% \frac{2^{1-2d_2}\Gamma(d_2)\Gamma((3/2)-d_2)}{\sqrt \pi\;(1-d_2)},\\%
%\frac{\Gamma(d_2)\Gamma(2-2d_2)}{\Gamma(2-d_2)}\\
%= \frac{C(d_1)(2-d_2)} {2C_*},\\
C^*&=&\int_0^1 x^{d_1-1}(1-x)^{1-2d_1}(1+x)\bigg\{\int_0^1
\frac{y^{d_2-1}(1-y)^{1-2d_2}(1+y)}{1+xy}\;\d y\bigg\}\d x.
% \frac{2^{1-2d_2}\Gamma(d_2)\Gamma((3/2)-d_2)}{\sqrt \pi \;(1-d_2)}
%\frac{\Gamma(d_2)\Gamma(2-2d_2)}{\Gamma(2-d_2)}\\
% &=&  ???  2^{1-2d_2} \pi^{-3/2} \sin(\pi d_1)
% \Gamma\Big(\frac{3}{2}-d_2\Big)
% \Gamma^{-1}\Big(\frac{2-d_2}{2}\Big).
\end{eqnarray*}
The variance of the noise is
\beq\label{e:sss}
 \sigma^2_\varepsilon = (2\pi)^{-1} C^* C(d_1) C(d_2)
 \sigma^2_{1,\varepsilon}\sigma^2_{2,\varepsilon}=
 \frac{\sin (\pi d_1)\sin (\pi d_2) C^*}{2\pi^3}.
\eeq
\end{proposition}

\noindent {\sc Proof.} \refeq{mix} implies that
\beq\label{e:bb}
 \varphi(x;d_1,d_2)= \frac{1}{C_*}\; (C(d_2)\varphi_1(x)
 F(-x;d_2)+C(d_1)\varphi_2(x) F(x;d_1)),
\eeq
 where
\beqn
 F(x;d): = \int_0^1 \frac{y^{d-1}
 (1-y)^{1-2d}(1+y)}{(1-xy)(1-y/x)}\; \d y,
\eeqn
\beqn
 C_* = \int_0^1\varphi_1(x)\bigg(\int_{-1}^0
 \frac{\varphi_2(y)}{1-xy}\;\d y\bigg)\d x =C(d_1)C(d_2)C^*.
\eeqn
% Here
%\begin{eqnarray*}
% \int_{-1}^0\frac{\varphi_2(y)}{1-xy}\;\d y
%&=&C(d_2)\int_0^1 y^{d_2-1} (1-y)^{1-2d_2} (1+xy)^{-1}(1+y)\d y\\
% &=& C(d_2)\bigg(\int_0^1 y^{d_2-1} (1-y)^{1-2d_2} (1+xy)^{-1}\d y\\
% &&+\int_0^1 y^{d_2} (1-y)^{1-2d_2} (1+xy)^{-1}\d y\bigg) \\
%&=&C(d_2)\bigg(\frac{\Gamma(d_2)\Gamma(2-2d_2)}{\Gamma(2-d_2)}\;
%F(1,d_2,2-d_2;-x)\\
% &&+\frac{\Gamma(d_2+1)\Gamma(2-2d_2)}{\Gamma(3-d_2)}\;
%F(1,d_2+1,3-d_2;-x)\bigg)\\
%&=&  ????????
%\end{eqnarray*}

Using equality
\beqn
 \frac{1+y}{(1-xy)(1-y/x)}=\frac{1}{(1-x)(1-y/x)}-\frac{x}{(1-x)(1-xy)},
\eeqn
 we have
\begin{eqnarray}
 F(x;d)&=& \frac{1}{1-x}\int_0^1 \frac{y^{d-1}
 (1-y)^{1-2d}}{1-y/x}\; \d y-\frac{x}{1-x}\int_0^1 \frac{y^{d-1}
 (1-y)^{1-2d}}{1-xy}\; \d y\nonumber\\
 &=& \frac{\Gamma(d)\Gamma(2-2d)}{\Gamma(2-d)}\frac{1}{1-x}\
 (F(1,d,2-d;1/x)-xF(1,d,2-d;x))\nonumber\\
 &=&\frac{\Gamma(d)\Gamma(2-2d)}{\Gamma(2-d)}\frac{G(x;d)}{1-x}.\label{e:cc}
%&=&\frac{2^{1-2d}\Gamma(d)\Gamma((3/2)-d)}{\sqrt \pi \;(1-d)} \;\frac{G(x;d)}{1-x} .\label{e:cc}
\end{eqnarray}
Now, \refeq{phif} follows from \refeq{bb} and \refeq{cc},
 whereas \refeq{sss} follows from \refeq{sigma2epsilon},
 \refeq{noise:fd}, \refeq{snoise:fd}.

% and equality
%\beqn
% C(d)\;\frac{2^{1-2d}\Gamma(d)\Gamma((3/2)-d)}{\sqrt \pi \;(1-d)}
% = \frac{2-d}{2}.
%\eeqn

 To finish the proof note that all the hypergeometric functions
appearing in the form of the mixture density are correctly
defined. \hfill $\Box$
\medskip

In the next proposition we present the asymptotics of
$\varphi(x;d_1,d_2)$ in the neighborhoods of 0 and $\pm 1$.

\begin{proposition} \label{p:FF}
Let the mixture density $\varphi$ be given in \refeq{phif}. Then
\begin{eqnarray}
\varphi(x;d_1,d_2)&\sim&
\begin{cases}
 \displaystyle\frac{\pi}{C^* \sin(\pi d_2)}\; x^{d_1+d_2-1},& \ \ x \to 0+,\\
 \displaystyle\frac{\pi}{C^* \sin(\pi d_1)}\; |x|^{d_1+d_2-1},& \ \ x \to 0-,
\end{cases}\label{e:phi11}\\
\varphi(x;d_1,d_2)&\sim&
\begin{cases}
 \displaystyle\frac{2^{1-2d_2}\pi}{C^* \sin(\pi d_2)}\;
 (1-x)^{1-2d_1},& \ \ x \to 1-,\\
 \displaystyle\frac{2^{1-2d_1}\pi}{C^* \sin(\pi d_1)}\; (1-x)^{1-2d_2},& \ \ x \to
 -1+.
\end{cases}\label{e:phi22}
\end{eqnarray}
\end{proposition}

\noindent {\sc Proof.} Applying identities (see Abramowitz and Stegun (1972))
\beqn
 F(a,b;c;1/x)= \Big(\frac{x}{x-1}\Big)^b
 F(b,c-a;c;\frac{1}{1-x}),\ \
F(a,b;c;1)=\frac{\Gamma(c)\Gamma(c-a-b)}{\Gamma(c-a)\Gamma(c-b)},
\eeqn
 we have that for $x\to 0+$
\begin{eqnarray*}
 G(-x,d_2)&=&F\Big(1,d_2,2-d_2;-\frac{1}{x}\Big)+x
 F(1,d_2,2-d_2;-x)\\
&=& \bigg(\frac{x}{1+x}\bigg)^{d_2} F(d_2,1-d_2,2-d_2; 1/(1+x))+x
 F(1,d_2,2-d_2;-x)\\
&\sim& x^{d_2} F(d_2,1-d_2,2-d_2;1)\\
&=& \frac{\Gamma(2-d_2)\Gamma(1-d_2)}{\Gamma(2-2d_2)} \ x^{d_2}\\
&=& \frac{\sqrt\pi\;\Gamma(2-d_2)}{2^{1-2d_2}\Gamma((3/2)-d_2)} \
x^{d_2},
\end{eqnarray*}
 and similarly for $x\to 0-$
\begin{eqnarray*}
 G(x,d_1)&\sim& \frac{\sqrt\pi\;\Gamma(2-d_1)}{2^{1-2d_1}\Gamma((3/2)-d_1)}\ |x|^{d_1}.
\end{eqnarray*}
% F\Big(1,d_1,2-d_1;\frac{1}{x}\Big)-x
% F(1,d_1,2-d_1;x)\\
%&=& \bigg(\frac{x}{x-1}\bigg)^{d_1} F(d_2,1-d_1,2-d_1; 1/(1-x))-x
% F(1,d_1,2-d_1;x)\\
%&\sim& |x|^{d_1} F(d_1,1-d_1,2-d_1;1)\\
%&=& \frac{\Gamma(2-d_1)\Gamma(1-d_1)}{\Gamma(2-2d_1)} \ |x|^{d_1}\\
%&=&....
%\end{eqnarray*}
This and equality \refeq{phif} imply
\begin{eqnarray*}
 \varphi(x;d_1,d_2)&\sim& C(d_1,d_2)\ \frac{\sqrt\pi\; \Gamma(2-d_2)}{2^{1-2d_2}\Gamma((3/2)-d_2)}\
 x^{d_1+d_2-1}\\
&=& \frac{C(d_1) C(d_2)} {C_*} \ \Gamma(d_2)\Gamma(1-d_2)  x^{d_1+d_2-1}\\
 &=& (C^*)^{-1} \frac{\pi}{\sin (\pi d_2)} \; x^{d_1+d_2-1},\ x\to 0+,\\
 \varphi(x;d_1,d_2)&\sim& %C(d_2,d_1)\
% \frac{\Gamma(2-d_1)\Gamma(1-d_1)}{\Gamma(2-2d_1)}\ x^{d_1+d_2-1}\\&=&
 (C^*)^{-1}  \frac{\pi}{\sin (\pi d_1)} \; |x|^{d_1+d_2-1},\ x\to
0-.
\end{eqnarray*}

For $x\to 1-$ we obtain
\begin{eqnarray*}
 \varphi(x;d_1,d_2)&\sim& 2C(d_1,d_2) F(1,d_2,2-d_2;-1)(1-x)^{1-2d_1}\\
 &=& C(d_1,d_2)\;\frac{\sqrt
 \pi\Gamma(2-d_2)}{\Gamma((3/2)-d_2)}\ (1-x)^{1-2d_1}\\
&=&  (C^*)^{-1} \frac{2^{1-2d_2}\pi}{\sin(\pi d_2)}\; (1-x)^{1-2d_1}
\end{eqnarray*}
 and similarly for $x\to -1+$
\begin{eqnarray*}
 \varphi(x;d_1,d_2)&\sim& 2C(d_2,d_1) F(1,d_1,2-d_1;-1)(1+x)^{1-2d_2}\\
 %&=& C(d_2,d_1)\;\frac{\sqrt
 %\pi\Gamma(2-d_1)}{\Gamma((3/2)-d_1)}\ (1+x)^{1-2d_2}\\
 &=& (C^*)^{-1} \frac{2^{1-2d_1}\pi}{\sin(\pi d_1)}\; (1-x)^{1-2d_2}.
\end{eqnarray*}

\begin{remark}
Clearly,
\begin{eqnarray*}
\varphi_1(x;d_1) &\sim&\begin{cases} C(d_1) x^{d_1-1},& \   x\to 0+,\\
 2 C(d_1) (1-x)^{1-2d_1},& \ x\to 1-,
\end{cases} \\
\varphi_2(x;d_2) &\sim&\begin{cases} C(d_2) |x|^{d_2-1},& \ x\to
0-,\\
 2C(d_2) (1+x)^{1-2d_2},&\ x\to -1+.
 \end{cases}
\end{eqnarray*}
Hence, by Proposition~\ref{p:FF}, the mixture density $\varphi$
associated with the product spectral density \refeq{ff} behaves as
$\varphi_1$ when $x$ approaches 1, and behaves as $\varphi_2$ when
$x$ approaches $-1$. % vanishes at the same rate, $(1-x)^{1-2d_1}$,
%as $\varphi_1$ when $x$ approaches 1 (and similarly at -1).
However, at zero, $\varphi$ behaves as $|x|^{d_1+d_2-1}$, i.e.\
both densities $\varphi_1$ and $\varphi_2$ count.
\end{remark}

Proposition~\ref{p:f12} allows us construct the mixture density
also in the case when the spectral density $f$ of aggregated
process has the form
\beqn
 f(\lambda)=\frac{1}{2\pi} \
  \Big(2\sin\frac{|\lambda|}{2}\Big)^{-2d} g(\lambda), \ \ 0<d<1/2,
\eeqn
 where %$C_*$ is a corresponding constant and
 $g(\lambda)$ is analytic spectral density on $[-\pi,\pi]$.
 In general, the existence of the mixture density associated with
 any analytic spectral density is not clear. For example, AR(1) is
 aggregated process only if the mixture density is the Dirac delta function,
 what is difficult to apply in practice. Similar inference concerns
 also the ARMA processes, i.e.\ rational spectral densities.
 Another class of spectral densities obtained by aggregating "non-degenerated"
 mixture densities is characterized in the following proposition.

\begin{proposition}\label{prop:analytic} A mixture density $\varphi_g$ is associated with
some analytic spectral density if and only if
%\beqn
% g(\lambda)=\frac{1}{2\pi}\int^1_0\frac{\varphi_g(x)}{|1-xe^{i\lambda}|^2}\; \d x.
%\eeqn
there exists $0<a_*<1$ such that ${\supp}(\varphi)\subset
[-a_*,a_*]$.
\end{proposition}

\noindent{\sc Proof.} For sufficiency, assume that there exists
$0<a_*<1$ such that ${\supp}(\varphi) \subset [-a_*,a_*]$. The
covariance function of the corresponding process satisfies
\begin{eqnarray*}
  |\gamma(h)| &\leq &\sigma^2_\varepsilon  \int_{-1}^1 \frac{|x|^{|h|}}{1-x^2}\; \varphi(x)\; \d x \\
  & =&\sigma^2_\varepsilon  \int_{-a_*}^{a_*} \frac{|x|^{|h|}}{1-x^2}\; \varphi(x)\; \d x  \\
 &\leq & \sigma^2_\varepsilon  a_*^{|h|}  \int_{-a_*}^{a_*} \frac{ \varphi(x)}{1-x^2}\; \d x \\
 &=& C a_*^{|h|},
\end{eqnarray*}
i.e.\  the covariance function decays exponentially to zero. This implies that the spectral density
$f(\lambda)=(2\pi)^{-1}\sum_{h=-\infty}^\infty \gamma(h)\e^{\i h\lambda}$ is analytic function on $[-\pi,\pi]$ (see, e.g., Bary (1964, p.\
80--82)).

To prove the necessity, assume that $f$ is an analytic function on
$[-\pi,\pi]$ or, equivalently, the corresponding covariance
function decays exponentially to zero, i.e.\ there exists $\theta
\in (0,1)$ and a constant $C>0$ such that $|\gamma(h)| \leq C
\theta^{|h|}$. Assume to the contrary that $\supp(\varphi) =
[-1,1]$.

Let $h\ge 0$ be an even integer. Then
$$
 |\gamma(h)|=\sigma^2_\varepsilon \int_{-1}^1 \frac{|x|^{h}}{1-x^2}\; \varphi(x) \; \d x \leq C \theta^h
$$
implies that
\begin{equation}
  \int_{-1}^1 \Big(\frac{|x|}{\theta}\Big)^h \frac{\varphi(x)}{1-x^2} \; \d x \leq
  C.
  \label{inter}
\end{equation}
Rewrite the last integral as
\begin{eqnarray*}
    \int_{-1/\theta}^{1/\theta} |x|^h \frac{ \varphi(\theta x)}{1-(\theta x)^2}\; \d x&=& \int_{-1}^1 |x|^h \frac{ \varphi(\theta x)}{1-(\theta x) ^2}
\; \d x  +
  \int_{-1/\theta}^{-1}  |x|^h \frac{ \varphi(\theta x)}{1-(\theta x) ^2} \; \d x\\
 && +\int_1^{1/\theta}  |x|^h \frac{ \varphi(\theta x)}{1-(\theta x)
^2} \; \d x =: I_1(h)+I_2(h)+I_3(h).
\end{eqnarray*}
For every $x\in[-1,1]$ we have
$$
 |x|^h \frac{\varphi(\theta x)}{1-(\theta x)^2}\leq
 \frac{\varphi(\theta x)}{1-(\theta x)^2}.
$$
Hence, by the dominated convergence theorem and \refeq{disp},
$I_1(h)\to 0$ as $h\to \infty$. The Fatou lemma, however, implies
that both integrals $I_2(h)$ and $I_3(h)$ tend to infinity as
$h\to \infty$:
$$
 \liminf_{h\to\infty} I_2(h) \geq  \int_1^{1/\theta}
 \frac{ \varphi(\theta x)}{1-(\theta x) ^2}  \liminf_{h\to\infty}  x^h  \; \d x  =\infty
$$
since $1/\theta>1$. This contradicts \eqref{inter}. \hfill$\Box$

\begin{example} Assume that the mixture density $\varphi_g$ has a uniform distribution on $[a,b]$, where
$-1<a<b<1$. By Proposition~\ref{prop:analytic}, the associated
spectral density is analytic function on $[-\pi,\pi]$ and can be
easily calculated:
\begin{eqnarray*}
 f_g(\lambda)&=&\frac{\sigma^2_\varepsilon}{2\pi(b-a)} \int_a^b \frac{\d x}{1-2x\cos
 \lambda+x^2}\\
&=&\frac{\sigma^2_\varepsilon}{2\pi(b-a)\sin|\lambda|}\ \bigg(\arctan\Big( \frac{b-\cos\lambda}{\sin|\lambda|}\Big)- \arctan\Big(
\frac{a-\cos\lambda}{\sin|\lambda|}\Big)\bigg), \ \lambda\ne 0,\pm \pi.
\end{eqnarray*}
$f_g(0)=\sigma^2_\varepsilon (2\pi)^{-1} (1-a)^{-1}(1-b)^{-1}, \ f_g(\pm\pi)=\sigma^2_\varepsilon (2\pi)^{-1} (1+a)^{-1}(1+b)^{-1}$.
%${\bf Linden (1999)???}
\end{example}

We  obtain the following corollary.
%{\bf [it is not very useful but it is then easier to formulate Prop. 4.2]}

\begin{corollary}\label{c:f12}
Let $\varphi_1(x;d)$ \refeq{phi1} and $\varphi_g(x)$ be the
mixture densities associated with spectral densities
$f_1(\lambda;d)$ \refeq{FARIMA} and analytic spectral density
$g(\lambda)$, respectively. Assume that ${\supp}(\varphi_g)
\subset [-a_*,0]$, $0<a_*<1$, and
% Let $X_t, t\in {\bf Z}$ be a process having spectral density of the form
\beq\label{e:fff12}
 f(\lambda) = \frac{1}{2\pi} \ \Big(2\sin\frac{|\lambda|}{2}\Big)^{-2d}
g(\lambda). %, \ \ where\ \  C_*:=\int_0^1\bigg(\int_{-a_*}^0
% \frac{\varphi_1(x;d)\varphi_g(y)}{1-xy}\;\d y\bigg)\d x.
\eeq
 Then the mixture density $\varphi(x)$, $x\in [-a_*,1]$ associated with $f$ is given by equality
\beqn%\label{e:mix}
 \varphi(x) =C_*^{-1}  \bigg(\varphi_1(x;d))
 \int_{-a_*}^{0}\frac{\varphi_g(y)}{(1-xy)(1-y/x)} \;\d y +
 \varphi_g(x) \int_{0}^{1}\frac{\varphi_1(y;d)}{(1-xy) (1-y/x)}  \;\d
 y\bigg),
\eeqn
 where
\beqn
  C_*:=\int_0^1\bigg(\int_{-a_*}^0
 \frac{\varphi_1(x;d)\varphi_g(y)}{1-xy}\;\d y\bigg)\d x.
\eeqn
\end{corollary}

%\begin{example} Assume that the mixture density {\bf [[Another example?]]}
%\end{example}

\section{The structure of aggregated process}\label{s:MA}

In order to make further inference about the aggregated process
$X_t$, e.g., estimation of the mixture density, limit theorems,
forecasting, etc., it is necessary to investigate more precise
structure of $X_t$. In particular, it is important to obtain the
linear (moving average) representation of the aggregated process.

\subsection{Behavior of spectral density of the aggregated process}

In this subsection we will study the behavior of the spectral
densities corresponding to the general class of semiparametric
mixture densities of the form (see Viano and Oppenheim (2004),
Leipus et al.\ (2006))
\beq\label{e:phi}
 \varphi(x) = (1-x)^{1-2d_1}(1+x)^{1-2d_2} \psi(x), \ \ 0<d_1<1/2,\
 0<d_2<1/2,
\eeq
 where $\psi(x)$ is continuous and nonvanishes at the points $x=\pm
 1$.
 As it is seen from the mixture densities appearing in Section~\ref{s:slm}, this form
 is natural, in particular \refeq{phi} covers
 the mixture density in \refeq{phif}. The
 corresponding spectral density behaves as a long memory
 spectral density.

% Throughout this section, without loss of generality, we assume the unit variances:
%\beqn

%\eeqn
%   {\bf [[we have to assume in this section the unit variances!]]}

%      First we formulate
% the lemma about the behaviour of the spectral density of the
% aggregated process $X_t$. Note, differently from Viano and
% Oppenheim (2004), we do not require the boundedness of function
% $\psi(x)$ on interval $[-1,1]$. In fact, $\psi(x)$ can have
% singularity points within $(-1,1)$ (see Example ...).

\begin{lemma} \label{lem:behav-spectr-dens}
 Let the density $\varphi(x)$ be given in \refeq{phi}, $\psi(x)$ is nonnegative
 function on $[-1,1]$ and continuous at
 the points $x=\pm 1$ with $\psi(\pm 1)\ne 0$. Then the following
 relations for the corresponding spectral density hold:
\begin{eqnarray}
 f(\lambda)&\sim& \frac{\sigma^2_\varepsilon\psi(1)}{2^{2d_2+1} \sin(\pi d_1)} \; |\lambda|^{-2d_1}, \quad |\lambda|\to 0,\label{e:f0}\\
 f(\lambda)&\sim& \frac{\sigma^2_\varepsilon\psi(-1)}{2^{2d_1+1}\sin(\pi d_2)}\; |\pi\mp\lambda|^{-2d_2},  \quad \lambda\to\pm \pi.\label{e:fp}
\end{eqnarray}
\end{lemma}

\noindent{\sc Proof.} Let $0<\lambda<\pi$. \refeq{phi},
\refeq{sp_density} and change of variables
$u=(x-\cos\lambda)/\sin\lambda$ lead to
\begin{eqnarray*}
 f(\lambda) &=&\frac{\sigma^2_\varepsilon}{2\pi}\int_{-1}^1 \frac{(1-x)^{1-2d_1}(1+x)^{1-2d_2}
 \psi(x)}{|1-x\e^{\i\lambda}|^2}
\; \d x \\
&=& \Big(2\sin\frac{\lambda}{2}\Big)^{-2d_1}
\Big(2\cos\frac{\lambda}{2}\Big)^{-2d_2} g^*(\lambda),
\end{eqnarray*}
 where \beqn
 g^*(\lambda)= \frac{\sigma^2_\varepsilon}{\pi}
\int_{-\cot\frac{\lambda}{2}}^{\tan\frac{\lambda}{2}} \frac{(\sin\frac{\lambda}{2}-u\cos\frac{\lambda}{2})^{1-2d_1}
(\cos\frac{\lambda}{2}+u\sin\frac{\lambda}{2})^{1-2d_2}}{1+u^2}\ \psi(u\sin\lambda+\cos\lambda)\; \d u. \eeqn By assumption of continuity at the
point 1, the function $\psi(u\sin\lambda+\cos\lambda)$ is bounded in some neighbourhood of zero, i.e.\ $\psi(u\sin\lambda+\cos\lambda)\le
C_1(\lambda_0)$ for $0<\lambda<\lambda_0$ and $\lambda_0>0$ sufficiently small. Hence,
$$
 \frac{(\sin\frac{\lambda}{2}-u\cos\frac{\lambda}{2})^{1-2d_1}
 (\cos\frac{\lambda}{2}+u\sin\frac{\lambda}{2})^{1-2d_2} }{1+u^2}\
 \psi(u\sin\lambda+\cos\lambda)\le \frac{C_2(\lambda_0)(1+|u|)^{1-2d_1}}{1+u^2}
$$
for $0<\lambda<\lambda_0$ and, by the dominated
convergence theorem, as $\lambda\to 0$,
\begin{eqnarray*}
g^*(\lambda)&\to& \frac{\sigma^2_\varepsilon\psi(1)}{\pi}\ \int_0^\infty \frac{u^{1-2d_1}}{1+u^2}\; \d u\\
&=& \frac{\sigma^2_\varepsilon\psi(1)}{2\pi}\ \Gamma(d_1) \Gamma(1-d_1)\\
&=&\frac{\sigma^2_\varepsilon\psi(1)}{2\sin(\pi d_1)}
\end{eqnarray*}
implying \refeq{f0}.
The same argument leads to relation \refeq{fp}.
 \hfill$\Box$

\begin{remark} Note, differently from Viano and Oppenheim (2004),
 we do not require the boundedness of function
 $\psi(x)$ on interval $[-1,1]$. In fact, $\psi(x)$ can have
 singularity points within $(-1,1)$, see \refeq{phi1}, \refeq{phi2}.
\end{remark}

\subsection{Moving average
representation of the aggregated process}

Any aggregated process admits an absolutely continuous
 spectral measure. If, in addition, its spectral density, say, $f(\lambda)$
 satisfies
\beq\label{e:reg1}
 \int_{-\pi}^\pi\log f(\lambda)\d \lambda >-\infty,
\eeq
 then the function
\beqn
 h(z)=\exp\Big\{\frac{1}{4\pi} \int_{-\pi}^\pi\frac{\e^{\i\lambda}+z}{\e^{\i\lambda}-z}\; \log f(\lambda)\d \lambda\Big\},
 \ \ |z|<1,
\eeqn
 is an outer function from the Hardy space $H^2$, does not vanish for $|z|<1$ and
 $f(\lambda)=|h(\e^{\i\lambda})|^2$. Then, by the Wold decomposition theorem,
 corresponding process $X_t$ is purely nondeterministic
 and has the MA($\infty$) representation (see Anderson (1971, Ch.\ 7.6.3))
\beq \label{e:MArepr}
 X_t=\sum_{j=0}^\infty \psi_j Z_{t-j},
\eeq
 where the coefficients $\psi_j$ are defined from the expansion of normalized outer function
 $h(z)/h(0)$, $\sum_{j=0}^\infty\psi_j^2<\infty$, $\psi_0=1$, and
 $Z_t=X_t-\widehat X_t$, $t=0,1,\dots$ ($\widehat X_t$ is the optimal linear
 predictor of $X_t$) is the innovation process, which is zero mean, uncorrelated,
 with variance
\beq \label{e:sigma22}
 \sigma^2=2\pi\exp\Big\{\frac{1}{2\pi}\int_{-\pi}^\pi\log f(\lambda)\d\lambda
 \Big\}.
\eeq
 By construction, the aggregated processes are Gaussian,
 implying that the innovations $Z_t$ are i.i.d.\ N$(0,\sigma^2)$ random variables.
\medskip

We obtain the following results.
% on the moving average representation of aggregated processes.

 \begin{proposition} Let the mixture density $\varphi$ satisfies the
 assumptions of Lemma~\ref{lem:behav-spectr-dens}. Then the
 aggregated process admits a moving average representation \refeq{MArepr},
%\beq\label{e:wold1}
% X_t=\sum_{j=0}^{\infty}\psi_j Z_{t-j},
%\eeq
 where the $Z_t$ are Gaussian i.i.d.\ random variables with zero mean and
variance \refeq{sigma22}.
%\beqn %\label{e:sigma2}
% \sigma^2=2\pi\exp\Big\{\frac{1}{2\pi}\int_{-\pi}^\pi\log f(\lambda)\d\lambda
% \Big\}
%\eeqn
% and the $\psi_j$ satisfy $\sum_{j=0}^\infty \psi_j^2<\infty$, $\psi_0=1$.
\end{proposition}

%{\bf [is it possible to say something about $\psi_j\sim$ in this case?]}

\noindent{\sc Proof.} We have to verify that \refeq{reg1} holds.
%It is well known that existence of the moving average representation
%\refeq{wold1} is guaranteed by the condition \beq\label{e:reg1}
% \int_{0}^\pi \log f(\lambda) \d \lambda > -\infty.
%\eeq
According to \refeq{f0}, $\log f(\lambda) \sim -C_1
 \log|1-\e^{\i\lambda}|$, $|\lambda|\to 0$, where $C_1>0$.
For any $\epsilon>0$ choose $0<\lambda_0\le \pi/3$, such that
\beqn
 - \frac{\log f(\lambda)}{C_1 \log|1-\e^{\i\lambda}|}-1\ge -\epsilon, \ \
 0<\lambda\le \lambda_0.
\eeqn
 Since $-\log|1-\e^{\i\lambda}|\ge 0$ for $0\le \lambda\le \pi/3$, we obtain
\beq\label{e:reg2}
 \int_0^{\lambda_0} \log f(\lambda)\d \lambda \ge (C_1-\epsilon)\int_0^{\lambda_0} \log |1-\e^{\i\lambda}|
 \d \lambda>-\infty
\eeq
 using the well known fact that $\int_{0}^\pi
\log|1-\e^{\i\lambda}| \d \lambda =0$. Using \refeq{fp} and the same argument, we get \beq\label{e:reg3}
   \int_{\pi-\lambda_0}^\pi \log f(\lambda) \d \lambda > -\infty.
\eeq

When $\lambda\in [\lambda_0, \pi-\lambda_0]$, there exist $0<L_1<L_2<\infty$ such that
$$
 L_1 \le \frac{1}{2\pi|1-x \e^{\i\lambda}|^2}\le L_2
$$
 uniformly in $x\in (-1,1)$. Thus, by \refeq{sp_density}, $L_1\le f(\lambda)\le L_2$
for any $\lambda\in [\lambda_0,\pi-\lambda_0]$, and therefore
\beq\label{e:reg4}
   \int_{\lambda_0}^{\pi-\lambda_0} \log f(\lambda) \d \lambda>-\infty.
\eeq
 \refeq{reg2}--\refeq{reg4} imply inequality \refeq{reg1}.
\hfill$\Box$

 \begin{lemma}\label{l:analytic}
If the spectral density $g$ of the aggregated process $X_t$ is
 analytic function, then $X_t$ admits representation
\beqn
 X_t=\sum_{j=0}^{\infty} g_j Z_{t-j},
\eeqn
 where the $Z_t$ are Gaussian i.i.d.\ random variables with zero mean and
variance
\beq\label{e:sigma2}
 \sigma_g^2=2\pi\exp\Big\{\frac{1}{2\pi}\int_{-\pi}^\pi\log
 g(\lambda)\d\lambda
 \Big\}
\eeq
 and the $g_j$ satisfy $|\sum_{j=0}^\infty g_j|<\infty$, $g_0=1$.
\end{lemma}

\noindent{\sc Proof.} From Proposition~\ref{prop:analytic}, there exists $0<a_*<1$ such that
\beq\label{e:ss}
 g(\lambda) = \frac{\sigma^2_\varepsilon}{2\pi}\int_{-a_*}^{a_*}
 \frac{\varphi_g(x)}{|1-x \e^{\i\lambda}|^2}\;  \d x.
\eeq
 For all $x\in[-a_*,a_*]$ and $\lambda\in[0,\pi]$ we have
$$
 \frac{1}{|1-x \e^{\i\lambda}|^2 } \geq C > 0.
$$
 This and \refeq{ss} imply $\int_{0}^\pi
\log g(\lambda) \d\lambda>-\infty$. Finally,
 $|\sum_{j=0}^\infty g_j|<\infty$ follows from representation
\beqn
 g(\lambda)= \frac{\sigma^2_g}{2\pi} \Big|\sum_{j=0}^\infty g_j \e^{\i j\lambda}\Big|^2
\eeqn
 and the assumption of analyticity of $g$.
\hfill$\Box$

\begin{proposition}\label{l:wold} Let $X_t$ be an aggregated process with spectral
density
\beqn
 f(\lambda)=\frac{1}{2\pi} \
  \Big(2\sin\frac{|\lambda|}{2}\Big)^{-2d} g(\lambda), \ \ 0<d<1/2,
\eeqn
 satisfying the assumptions of Corollary~\ref{c:f12}.
%\beqn
% f(\lambda)=\frac{1}{C_*} \Big(2\sin\frac{\lambda}{2}\Big)^{-2d}
% g(\lambda), \ \ 0<d<1/2,
%\eeqn
% where $g$ is analytic spectral density of some aggregated process.
 Then $X_t$ admits a representation \refeq{MArepr},
%\beq\label{e:ma}
% X_t=\sum_{j=0}^{\infty}\psi_j Z_{t-j},
%\eeq
 where the $Z_t$ are Gaussian i.i.d.\ random variables with zero mean and variance
\beqn%\label{e:sigma2}
 \sigma^2=2\pi\exp\Big\{\frac{1}{2\pi}\int_{-\pi}^\pi\log f(\lambda)\d\lambda\Big\}=
 \exp\Big\{\frac{1}{2\pi}\int_{-\pi}^\pi\log g(\lambda)\d\lambda\Big\} =\frac{\sigma^2_g}{2\pi}
\eeqn and the $\psi_j$ satisfy %$\psi_0=1$,
\beq\label{e:phi_as}
 \psi_j \sim \frac{\sum_{j=0}^\infty g_j}{\Gamma(d)} \; j^{d-1}, \ \ \psi_0=1.
\eeq
 Here, the $g_j$ are given in Lemma \ref{l:analytic}.
\end{proposition}

\noindent{\sc Proof.} We have
 \beqn
 \frac{1}{2\pi}\Big(2\sin\frac{|\lambda|}{2}\Big)^{-2d}=
  \frac{1}{2\pi}\bigg|\sum_{j=0}^\infty h_j \e^{\i j\lambda}\bigg|^2\
 \ {\rm with} \ \ h_j =\frac{\Gamma(j+d)}{\Gamma(j+1)\Gamma(d)}
\eeqn
 and, recall,
\beqn
 g(\lambda)=\frac{\sigma^2_g}{2\pi}\bigg|\sum_{j=0}^\infty g_j \e^{\i j\lambda}\bigg|^2,\
  \ \sum_{j=0}^\infty g^2_j <\infty
\eeqn
 since, by Lemma~\ref{l:analytic}, $\int_{-\pi}^\pi \log
 g(\lambda)\d\lambda>-\infty$.
 On the other hand, $\int_{-\pi}^\pi \log f(\lambda)\d\lambda>-\infty$ implies
\beqn
 f(\lambda)=\frac{1}{2\pi}\Big|\sum_{j=0}^\infty \tilde\psi_j \e^{\i j\lambda}\Big|^2, \ \
 \sum_{j=0}^\infty\tilde \psi_j^2<\infty
\eeqn
 and, by uniqueness of the representation,
\beqn
 \tilde\psi_k= \frac{\sigma_g}{\sqrt{2\pi}}\sum_{j=0}^k h_{k-j} g_j.
\eeqn
 It easy to see that,
\beq\label{e:dct}
  \sum_{j=0}^k h_{k-j} g_j\sim h_k \sum_{j=0}^\infty g_j\sim C_2 k^{d-1},
\eeq
 where $C_2= \Gamma^{-1}(d) \sum_{j=0}^\infty g_j$. Indeed, taking into
 account that $h_k\sim \Gamma^{-1}(d) k^{d-1}$, we can write
\beqn
 \sum_{j=0}^k h_{k-j} g_j = \Gamma^{-1}(d)k^{d-1}\sum_{j=0}^\infty a_{k,j}
 g_j,
\eeqn
 where $a_{k,j}= h_{k-j} \Gamma(d) k^{1-d} {\bf 1}_{\{j\le k\}}\to
 1$ as $k\to\infty$ for each $j$. On the other hand, we have $|a_{k,j}|\le C
 (1+j)^{1-d}$ uniformly in $k$ and, since the $g_j$ decay exponentially fast,
 the sum $\sum_{j=0}^\infty (1+j)^{1-d} |g_j|$ converges and
 the dominated convergence theorem applies to obtain \refeq{dct}.

 Hence, we can write
\beqn
 f(\lambda)=\frac{\sigma^2_g}{(2\pi)^2}\Big|\sum_{j=0}^\infty \psi_j \e^{\i j\lambda}\Big|^2, \
 \ \psi_0=1,
\eeqn
 where $\psi_j=\tilde \psi_j\sqrt {2\pi}/\sigma_g\sim C_2 j^{d-1}$ and \refeq{MArepr} follows.
% \beqn
% X_t=\sum_{j=0}^\infty\tilde \psi_j \tilde Z_{t-j},
%\eeqn
% where $\tilde Z_t\sim {\rm i.i.d.}\; {\rm N}(0,1)$ and the proposition follows noting that
%\beqn
% \tilde \psi_0=\frac{\sigma_g}{\sqrt{2\pi}}=
% \exp\Big\{\frac{1}{4\pi}\int_{-\pi}^\pi\log g(\lambda)\d\lambda\Big\}.%=\frac{\sigma_g}{\sqrt{2\pi}}.
%\eeqn
 \hfill$\Box$

\section{Appendix. Mixture density associated with FI($d$) spectral density}

\begin{proposition}\label{p:DCO}
 Mixture density associated with FI($d$) spectral density
\beqn
 f(\lambda;d)=\frac{1}{2\pi}\; \Big(2\sin\frac{|\lambda|}{2}\Big)^{-2d},
 \ 0<d<1/2,
\eeqn
 is given by equality
\beq\label{e:dco1}
 \varphi(x)=
 C(d) x^{d-1} (1-x)^{1-2d} (1+x) {\bf 1}_{[0,1]}(x),
\eeq
 where
\beqn
 C(d)=\frac{\Gamma(3-d)}{2\Gamma(d)\Gamma(2-2d)}= 2^{2d-2}\frac{\sin(\pi d)}{\sqrt \pi} \frac{\Gamma(3-d)}{\Gamma ((3/2)-d)}.
\eeqn
 The variance of the noise is
\beqn
  \sigma^2_\varepsilon= \frac{\sin(\pi d)} {C(d)\pi}.
\eeqn
 \end{proposition}

\noindent {\sc Proof.} Equality
$$
f(\lambda;d)=2^{-d-1}\pi^{-1} (1-\cos \lambda)^{-d}
$$
 implies that $1-\cos \lambda = (\pi
2^{d+1}f(\lambda;d))^{-1/d}$. Hence, rewriting
\begin{eqnarray*}
 |1-x\e^{\i\lambda}|^2=(1-x)^2 \Big(1+\frac{2x}{(1-x)^2}\; (1- \cos
\lambda)\Big),
\end{eqnarray*}
 and assuming ${\supp}(\varphi)=[0,1]$, we obtain that the spectral density
 of aggregated process is of the form
\begin{eqnarray}
\frac{\sigma^2_\varepsilon}{2\pi} \int_0^1
 \frac{\varphi(x)}{|1-x\e^{\i\lambda}|^2}\; \d x
 &=& \frac{\sigma^2_\varepsilon}{2\pi}
 \int_{0}^1 \frac{\varphi(x)}
 {(1-x)^2 \big(1+\frac{2x}{(1-x)^2} (1- \cos \lambda)\big)} \;\d x \nonumber\\
 &=& \frac{\sigma^2_\varepsilon}{2\pi} \int_{0}^1
 \frac{\varphi(x)}{(1-x)^2\big(1+\frac{2x}{(1-x)^2}(\pi
 2^{d+1}f(\lambda;d))^{-1/d}\big)} \;\d x \label{e:olia}.
\end{eqnarray}
 The change of variables $y^{1/d}=2x/(1-x)^2$ implies
\beq\label{e:galia1}
 \d y =
 \frac{d 2^d x^{d-1}(1+x)}{(1-x)^{2d+1}}\; \d x.
\eeq

Consider the density $\varphi$ defined by
\beq\label{e:galia2}
 \d y =\frac{\varphi(x)}{\tilde C(d) (1-x)^2}\; \d x,
\eeq
 where $\tilde C(d)$ is some constant. Then \refeq{olia}
 becomes
\begin{eqnarray}
 \frac{\sigma^2_\varepsilon}{2\pi} \int_0^1
 \frac{\varphi(x)}{|1-x\e^{\i\lambda}|^2}\; \d x&=&
 \frac{\sigma^2_\varepsilon\tilde C(d)}{2\pi} \int_0^{\infty} \frac{\d y}
 {(1+y^{1/d} (\pi 2^{d+1}f(\lambda;d))^{-1/d})}\nonumber\\
 &=& f(\lambda;d) \sigma^2_\varepsilon 2^d \tilde C(d) \int_0^\infty \frac{\d
 z}{1+z^{1/d}}\nonumber\label{e:ba}
\end{eqnarray}
 after the change of variables $z=\frac{y}{\pi 2^{d+1}f(\lambda;d)}$.
Therefore, FI($d$) is an aggregated process and, by \refeq{galia1}--\refeq{galia2},
 the mixture density has a form
\beq\label{e:md}
 \varphi(x) = \tilde C(d)d 2^d x^{d-1}(1-x)^{1-2d} (1+x),
\eeq
 and the variance of the noise is (see formula 6.1.17 in Abramowitz and Stegun (1972))
\begin{eqnarray*}
 \sigma^2_\varepsilon&=&2^{-d} (\tilde C(d))^{-1}\bigg(\int_0^\infty \frac{\d
 z}{1+z^{1/d}}\bigg)^{-1}\\
 &=&2^{-d} (d \tilde C(d))^{-1}(B(d,1-d))^{-1}\\
 &=&2^{-d} (d \tilde C(d))^{-1}\; \frac{\sin (\pi d)}{\pi}.
\end{eqnarray*}

Finally, it remains to calculate the constant $\tilde C(d)$ to ensure that the mixture density
 $\varphi$ given in \refeq{md} integrates to one over the interval $[0,1]$. We have
\begin{eqnarray*}
\int_0^1 \varphi (x) \d x
&=& \tilde C(d) d 2^d \Big(\int_0^1 x^{d-1}(1-x)^{1-2d} \d x + \int_0^1 x^d (1-x)^{1-2d} \d x\Big)  \\
&=& \tilde C(d) d 2^d \Big(B(d,2-2d)+B(d+1,2-2d)\Big)\\
%&=& \tilde C(d) d 2^{1+d} (2-d)^{-1} B(d,2-2d) \\
&=& \tilde C(d) d 2^{1+d}\; \frac{\Gamma(d)\Gamma (2-2d))}{(2-d)\Gamma(2-d)} \\
&=& \tilde C(d) d 2^{2-d} \frac{\sqrt\pi}{\sin(\pi d)}
\frac{\Gamma((3/2)-d)}{\Gamma(3-d)}.
\end{eqnarray*}
 Hence,
\beqn
 \tilde C(d)=\frac{1}{d2^{d+1}} \frac{\Gamma(3-d)}{\Gamma(d)\Gamma(2-2d)}= \frac{2^{d-2}}{d}\frac{\sin(\pi d)}{\sqrt \pi} \frac{\Gamma(3-d)}{\Gamma ((3/2)-d)}
\eeqn
 and $C(d)=\tilde C(d) d 2^d$.
\hfill $\Box$
\medskip

\newpage

\section*{References}

\bigskip

\begin{description}

\itemsep -.04cm

\item Abramowitz, M., Stegun, I.A. (1972) {\em Handbook of
Mathematical Functions}. Dover, New York.

\item Anderson, T.W. (1971) {\em The Statistical Analysis Of Time Series}. John Wiley \& Sons, New York.

\item Bary, N.K. (1964) {\em A Treatise On Trigonometric Series,
Vol.\ I}. Pergamon Press, Oxford.

%\item Bhansali, R.J., Giraitis, L. and Kokoszka, P.S. (2005) Approximations
%and limit theory for quadratic forms of linear processes. Preprint.

\item Chong, T.T. (2006) The polynomial aggregated AR(1) model.
{\em Econometrics Journal} {\bf 9}, 98--122.

\item Dacunha-Castelle, D. and Oppenheim, G. (2001) Mixtures, aggregations
and long memory. Universit\' e de Paris-Sud, Mathematiques. Preprint 2001-72.

\item Dacunha-Castelle, D. and Fermin, L. (2006) Disaggregation of
long memory processes on $C^\infty$ class. {\em Electronic
Communications in Probability} {\bf 11}, 35--44.

\item Haubrich, J.G. and Lo, A.W. (2001) The sources and nature of
long-term dependence in the business cycle. {\em Econom. Rev.}
(Fed. Reserve Bank Cleveland, Q II) {\bf 37}, 15--30.

\item Leipus, R., Oppenheim, G., Philippe, A. and Viano, M.-C. (2006)
 Orthogonal series density estimation in a disaggregation scheme. {\em
 Journal of Statistical Planning and Inference} {\bf 136}, 2547--2571.

\item Oppenheim, G. and Viano, M.-C. (2004)
Aggregation of random parameters Ornstein-Uhlenbeck or AR processes: some
convergence results. {\em Journal of Time Series Analysis} {\bf 25}, 335--350.

\item{Szeg\"o, G.} (1959) {\em Orthogonal Polynomials}.
American Mathematical Society, New York.

\item{Tolstov, G.P.} (1976) {\em Fourier Series}.
Dover, New York.

\item{Zaffaroni, P.} (2004) Contemporaneous aggregation of linear
dynamic models in large economies. {Journal of Econometrics} {\bf
 120}, 75–-102.

\end{description}

\end{document}